\documentclass{amsart}

\usepackage{amsmath,amsthm,amssymb,amscd,mathabx,accents,enumitem,mathtools}
\usepackage{amsfonts}
\usepackage{rotating}
\usepackage{euscript}
\usepackage{pst-node}
\usepackage{epsfig,verbatim}
\usepackage{tikz-cd}

\def\be{\begin{equation}}
\def\ee{\end{equation}}

\def\C{{\mathbb C}} 
\def\f{\EuScript}
\def\N{{\mathbb N}} 
\def\P{{\mathbb P}}

\def\Q{{\mathbb Q}}

\def\phi{{\varphi}}

\def\Aut{{\rm Aut}}

\def\LCM{{\rm LCM }}

\def\bp{\begin{proposition}}
\def\ep{\end{proposition}}

\def\bt{\begin{theorem}}
\def\et{\end{theorem}}
\def\bpr{\begin{problem}}
\def\epr{\end{problem}}
\def\bco{\begin{conjecture}}
\def\eco{\end{conjecture}}
\def\be{\begin{equation}}
\def\bee{\begin{equation*}}
\def\l{\label}

\def\ee{\end{equation}}
\def\eee{\end{equation*}}
\def\bl{\begin{lemma}}
\def\el{\end{lemma}}
\def\bc{\begin{corollary}}
\def\ec{\end{corollary}}

\def\bd{\begin{definition}}
\def\ed{\end{definition}}

\def\hat{\widehat}

\mathtoolsset{showonlyrefs}

\newtheorem{theorem}{Theorem}[section]
\newtheorem{lemma}[theorem]{Lemma}
\newtheorem{corollary}[theorem]{Corollary}
\newtheorem{proposition}[theorem]{Proposition}
\newtheorem{problem}[theorem]{Problem}
\newtheorem{conjecture}[theorem]{Conjecture}

\begin{document}

\begin{abstract}
We formulate some problems and conjectures about semigroups of rational functions under composition.  
The considered problems arise in different contexts, but most of  them are united by a certain relationship
  to the concept of  amenability.  

\end{abstract}

\title[Semigroups of rational functions under compositions]{Semigroups of rational functions: 
 some  problems and conjectures}
\author{Fedor Pakovich}
\thanks{
This research was supported by ISF Grant No. 1432/18}
\address{Department of Mathematics, Ben Gurion University of the Negev, Israel}
\email{
pakovich@math.bgu.ac.il}

\maketitle

\section{Introduction}
The goal of this paper is twofold. First, we propose a number of problems and conjectures about semigroups of rational functions under composition. Second, and not less importantly, we create bridges  between some of these problems, which arise in different contexts and apparently are not related. 
 Most of the problems discussed in the paper can be linked  in some way to the concept of amenability. Nevertheless,  sometimes the corresponding link is rather indirect. The paper does not pretend to be a review of existing papers about semigroups of rational functions,  and the choice of considered problems is  determined  only by the tastes and interests of the author.

We recall that a semigroup $S$ is called {\it left amenable} if it admits a finitely additive probability
measure $\mu$, defined on all the subsets of $S$, which is left invariant in the following sense. For all $\f T\subseteq S$ and $A \in S$ the equality 
\be \l{perd} \mu(A^{-1}\f T) = \mu(\f T)\ee holds, where the set $A^{-1} \f T$ is defined by the formula  $$A^{-1} \f T=\{W \in S\, | \,AW \in \f T\}.$$ 
Equivalently, $S$ is left amenable if there is a mean on $l_{\infty}(S)$, which is invariant under the natural
left action of $S$ on the dual space $l_{\infty}(S)^*$ (see e.g. \cite{pater}). The right amenability is defined similarly. A semigroup is called {\it amenable} if 
there exists a mean on $l_{\infty}(S)$, which is invariant under the left and the right action of $S$ on $l_{\infty}(S)^*$. By the theorem of Day (see \cite{day0}, \cite{day}), this is equivalent to the condition that $S$ is left and right amenable,  and in this paper we will use the last condition as the definition of amenability.

To our best knowledge, amenability in semigroups of rational functions started to be studied only recently in the  papers \cite{peter}, \cite{peter2}, \cite{amen}. Some close questions were considered also in the paper \cite{bell}. 
Let $P$ be a rational function. The results of the above mentioned papers  show that the most interesting semigroups of rational functions related to the concept of amenability are 
semigroups $C(P),$ $C_{\infty}(P)$, and $E(P)$, which can be defined correspondingly as semigroups consisting of rational functions commuting with $P$, commuting with some {\it iterate} of $P$, and sharing a measure of maximal entropy with $P$.  On the other hand,  the most interesting property of amenable semigroups in the context of semigroups of rational functions is, probably, the property of being reversible.  
We recall that a semigroup $S$ is called {\it left reversible} 
if for all $A,B\in S$ there exist $X,Y\in S$ such that \be \l{reva} AX=BY.\ee The right reversibility is defined similarly. It is well-known and follows easily from the definition  that any left (resp. right) amenable semigroup is left (resp. right) reversible (see e.g. \cite{pater}).

What is said above suggests to consider problems related to amenable semigroups of rational functions in conjunction  with   problems related to reversible semigroups and semigroups $C(P),$ $C_{\infty}(P)$, $E(P)$. This paper is organized accordingly and has the following structure. In the second section, after recalling the classical theorem of Ritt  about commuting rational functions, we review recent results about $C(P)$ and propose some related problems. In the third section,  we discuss some conjectures about semigroups  $C_{\infty}(P)$ and $E(P)$. In the fourth section, we discuss problems and conjectures concerning reversible semigroups of rational functions. In particular, we relate the left reversibility with some problems arising in arithmetical dynamics. 
Finally, in the fifth section we propose some conjectures describing amenable semigroups of rational functions, generalizing previous results obtained in the polynomial case.

\section{Semigroups $C(P)$}

 For a rational function $P$ of degree at least two, we denote by $C(P)$  
the collection of rational functions, including rational functions of degree one, commuting with $P$. It is clear that $C(P)$ is  a semigroup. The subsemigroup of $C(P)$ consisting of   M\"obius
transformations will be denoted by $\Aut(P).$  It is easy to see that $\Aut(P)$ is a group. Moreover, since elements of $\Aut(P)$ permute fixed points of $P^{\circ k},$ $k\geq 1$, and any M\"obius 
transformation is defined by its values at any three points, the group $\Aut(P)$ is finite. We will call a rational function  {\it special} if it is 
either a Latt\`es map or it is conjugate to  $z^{\pm n},$ $n\geq 2$, or $\pm T_n,$ $n\geq 2$, where $T_n$ is the $n$th Chebyshev polynomial. 

The central result about commuting rational functions is the theorem of Ritt  (see \cite{r} and also \cite{e2}, \cite{fin}, \cite{rev}). In brief, it states that  if rational functions $X$ and $P$ of degree at least two commute, then either they both are special or they have an iterate in common. Moreover, 
the Ritt theorem provides a full description of pairs of commuting special functions, implying a description of the semigroup $C(P)$ in case  $P$ is special. Nevertheless, the only information about the structure of $C(P)$ for  a non-special rational function $P$ provided by the Ritt theorem is that every $X\in C(P)$ has a common iterate with $P$. Further results about $C(P)$ were obtained in the recent papers \cite{fin}, \cite{rev}. Below, we shortly describe these results and formulate some related problems.

The following theorem was obtained in the paper \cite{fin} as a corollary of a more general theorem concerning semiconjugate rational functions. It can be regarded as a ``semigroup'' counterpart of the Ritt theorem and implies the latter theorem in its part concerning non-special rational functions.   

\pagebreak

\bt \l{gr0} Let $P$ be a non-special rational function of degree at least two. Then there 
exist {\it finitely many} rational functions $X_1, X_2,\dots, X_r$ such that a rational function $X$ belongs to $C(P)$
if and only if 
\be \l{rep} X=X_j\circ P^{\circ k}\ee
for some  $j,$ $1\leq j \leq r,$  and  $k\geq 0.$
\et  

To see that Theorem \ref{gr0}  implies the Ritt theorem, let us observe that  if $X$ commutes with  $P$, then 
any iterate  $X^{\circ l}$, $l\geq 1,$ does. Thus,  by the Dirichlet box principle, 
there exist  $l_2> l_1\geq 1$ such that 
$$X^{\circ l_1}=X_j\circ P^{\circ k_1}, \ \ \ \ \ \ \ X^{\circ l_2}=X_j\circ P^{\circ k_2}$$ for the same $j$ and 
some  $k_1>k_2\geq 1$, implying that 
\be \l{yuy} X^{\circ l_2}=X^{\circ l_1}\circ  P^{\circ (k_2-k_1)}.\ee
Since $X$ and $P$ commute, it follows from \eqref{yuy} that 
\be \l{yuy1} X^{\circ l_2}= P^{\circ (k_2-k_1)}\circ X^{\circ l_1},\ee
implying that 
\be \l{aob} X^{\circ (l_2-l_1)}=P^{\circ (k_2-k_1)}.\ee 
The passage from \eqref{yuy1} to \eqref{aob} is possible since the semigroup $\C(z)$ is obviously {\it right cancellative}, that is, the equality $$X\circ A=Y\circ A,$$ where $X,Y,A\in \C(z),$ implies that $X=Y.$ Notice, however,  that $\C(z)$ is not left cancellative.

It was shown in \cite{rev} that with the semigroup $C(P)$ one can associate a finite {\it group} as follows. For a non-special rational function of degree at least two $P$, we define an equivalence relation 
$\underset{P}{\sim}$ on the semigroup $C(P)$, setting  
$Q_1\underset{P}{\sim} Q_2$ if 
\be \l{zek} Q_1\circ P^{\circ l_1}=Q_2\circ P^{\circ l_2}\ee  
for some $l_1\geq 0,$ $l_2\geq 0$. It follows easily from the right cancellativity of  $\C(z)$ that if $\bf{X}$ is an equivalence class of $\underset{P}{\sim}$ and $X_0\in \bf{X}$ is a function of minimum possible degree, then 
every $X\in \bf{X}$ has the form $X=X_0\circ P^{\circ l},$ $l\geq 0.$ 
Moreover, the multiplication of classes induced by the functional composition of their representatives
provides $ C_B/\underset{B}{\sim}$ with the structure of a {\it finite group}.

In more detail, let us recall that a  {\it congruence} on a  semigroup is an equivalence relation that is compatible with the semigroup operation. In this notation, 
the following statement holds.

\bt \l{gr} Let $P$ be a non-special rational function of degree at least two. Then the relation $\underset{P}{\sim}$ is a congruence on the semigroup $C(P)$, and the quotient semigroup is a finite  group. 
\et

It is clear that describing the semigroup $C(P)$ is equivalent to describing the corresponding group, which will be denoted by $G_P$. On the other hand, one may expect that the analysis of $G_P$ may have some advantages in view of the presence of the group structure. Thus, probably, the most interesting problem concerning the semigroups $C(P)$ is following. 

\bpr \l{pr1} Which finite groups occur as groups $G_P$ for non-special rational functions $P$?
\epr

Note that the group $G_P$ is trivial if and only if any element of $C(P)$ is an iterate of $P.$ In particular, if $P=Q^{\circ l}$ for some $Q\in \C(z),$ then $G_P$ is non-trivial since $Q$ belongs to $C(P)$. The group $G_P$ is also non-trivial whenever the group $\Aut(P)$ is non-trivial, since $G_P$ contains an isomorphic copy of $\Aut(P)$ (see \cite{rev} for more detail).

Let us mention some known results related to Problem \ref{pr1}. First, if $P$ is a non-special {\it polynomial}, then $ C(P)=\langle \Aut(P),R\rangle$ for some  $R\in C(P)$.  Correspondingly, the group $G_P$  is metacyclic (see \cite{rev}, Section 6.2, and \cite{amen}, Section 7.1). 
Further, if  a non-special rational function $P$ is {\it indecomposable}, that is, cannot be represented in the form $P=V\circ U$, where $U$ and $V$ are rational functions of degree greater than one, then  
$G_P$  is isomorphic to $\Aut(P)$. 
Equivalently, $X\in C(P)$ if and only if $X=\mu\circ P^l$ for some $\mu\in \Aut(P)$ and $l\geq 1$ (see \cite{rev}, Section 6.1). Thus, for an indecomposable rational function $P$ the group $G_P$ is one of the five finite subgroups   of $\Aut(\C\P^1)$. Moreover, every finite subgroup of $\Aut(\C\P^1)$ can be realized as the group $\Aut(P)$ for some rational function $P$ (see \cite{dm}).

It was shown in \cite{rev} that calculating  $G_P$ reduces to calculating the generators of the fundamental group of a special graph  $\Gamma_P$ associated with $P$, which can be described as follows. Let $P$ be a rational function. A rational function $\hat P$ is called an {\it elementary transformation} of $P$
if there exist rational functions $U$ and $V$ such that  $P=V\circ U$ and $\hat P=U\circ V$. We say that rational functions 
 $P$ and $A$ are  {\it equivalent} and write $A\sim P$  if there exists 
a chain of elementary transformations between $P$ and $A$.
Since for any M\"obius transformation $\mu$ the equality
$$P=(P\circ \mu^{-1})\circ \mu$$ holds, 
the equivalence class $[P]$ of a rational function $P$ is a union of conjugacy classes.
Thus, the relation $\sim$ can be considered as a 
weaker form of the classical conjugacy relation. 
The graph $\Gamma_P $ is defined as a multigraph whose vertices
are in a one-to-one correspondence with some fixed representatives $P_i$ of conjugacy classes in $[P]$, and whose multiple edges connecting the vertices  corresponding to $P_i$ to $P_j$ 
are in a one-to-one correspondence with
solutions of the system  
$$P_i=V\circ U, \ \ \ P_j=U\circ V$$
in rational functions. Since the equivalence class $[P]$ 
contains infinitely many conjugacy classes if and only if 
$P$ is a flexible Latt\`es map (\cite{rec}), for any non-special rational function $P$ the graph  $\Gamma_P $ is finite. 

 It follows from the relation between $G_P$ and the fundamental group of  $\Gamma_P $  that the groups  $G_P$ and  $G_{P'}$ are isomorphic  whenever $P\sim P'$ (\cite{rev}). This fact permits to reveal   reasons for the non-triviality of $G_P$ by studying the graph $G_P$, which represents the totality of all functions from  class $[P]$ together with their decompositions and automorphisms. 
As an example of such an approach, we mention the following interesting fact. The group $\Gamma_P$ is always  non-trivial whenever for some $P'\sim P$  the group $\Aut(P')$ is non-trivial. Thus, if under these circumstances the group $\Aut(P)$ itself is trivial, the semigroup $C(P)$ necessarily contains functions of degree at least two that are not iterates of $P$  (see   \cite{rev} for more detail).

In view of the relation between $G_P$ and $\Gamma_P$, a ``combinatorial'' counterpart of Problem \ref{pr1} is the following problem.

\bpr \l{pr11} Which finite graphs occur as graphs $\Gamma_P$  for non-special rational functions $P$?
\epr

\section{Semigroups $C_{\infty}(P)$ and $E(P)$}

Let us define the sets $C_{\infty}(P)$  and $\Aut_{\infty}(P)$ by the formulas 
\be \l{defin} C_{\infty}(P)= \bigcup_{i=1}^{\infty}C(P^{\circ k}), \ \ \ \ \Aut_{\infty}(P)=\bigcup_{k=1}^{\infty} \Aut(P^{\circ k}).\ee
Since obviously  
\be \l{lcm} C(P^{\circ k}),\  C(P^{\circ l})\subseteq C(P^{\circ \LCM(k,l)})\ee 
and 
  \be \l{xre} \Aut(P^{\circ k}),\   \Aut(P^{\circ l})\subseteq  \Aut(P^{\circ \LCM(k,l)}),\ee 
the set $C_{\infty}(P)$ is a semigroup, and the set  $\Aut_{\infty}(P)$ is a group. Let us notice that 
a rational function $X$ of degree at least two belongs to  $C_{\infty}(P)$ if and only if $X$ and $P$ share an iterate. Indeed, the ``only if'' part follows from the Ritt theorem.  On the other hand, if there exist $k,l\in \N$ such that $X^{\circ k}=P^{\circ l},$ then $X$ obviously 
commutes with $P^{\circ l}.$

We conjecture that 
the following analogue of Theorem \ref{gr0} holds for semigroups  $C_{\infty}(P).$

\bco \l{co3}  Let $P$ be a non-special rational function of degree at least two. Then there 
exist {\it finitely many} rational functions $X_1, X_2,\dots, X_r$ such that  a rational function $X$ belongs to $C_{\infty}(P)$ if and only if 
\be \l{rep2} X=X_j\circ P^{\circ k}\ee
for some  $j,$ $1\leq j \leq r,$  and  $k\geq 0.$ 
\eco 

Conjecture \ref{co3} is equivalent to the following conjecture. 

\bco \l{co4}  Let $P$ be a non-special rational function of degree at least two. 
Then $C_{\infty}(P)=C(P^{\circ s})$ for some $s\geq 1.$ 
\eco

Indeed, applying Theorem \ref{gr0} to $P^{\circ s}$, we see  that Conjecture \ref{co4} implies Conjecture \ref{co3}. On the other hand, if   $X_1,X_2,\dots, X_r$ are rational functions satisfying the conclusion of 
Conjecture \ref{co3}, and  $X_i$, $1\leq i \leq r,$ commutes with $P^{\circ k_i}$, $k_i\geq 1,$ then $X_i$ also commutes with $P^{N},$ where $N=\LCM(k_1,k_2,\dots,k_r).$  
Thus, $C_{\infty}(P)\subseteq C(P^{\circ N}),$ implying that 
$C_{\infty}(P)=C(P^{\circ N}).$

Let us recall that by the results of the papers \cite{flm} and \cite{l}, 
 for every rational function $P$ of degree $n\geq 2$ there exists a unique probability measure $\mu_P$ on $\C\P^1$, which is invariant under $P$, has support equal to the Julia set $J_P$, and achieves maximal entropy $\log n$ among all $P$-invariant probability measures. 
For a  rational function  $P$ 
of degree at least two, we denote by $\mu_P$ the measure of 
 maximal entropy for $P$, and by $E(P)$ the set of rational functions $Q$ of degree at least two such that $\mu_Q=\mu_P$, completed by  $\mu_P$-invariant M\"obius transformations. The set $E(P)$ is a semigroup (see e.g. \cite{amen}).   

Algebraic conditions for non-special rational functions $X$ and $P$ to share a measure of maximal entropy were obtained in the papers \cite{lev}, \cite{lp}, and can be formulated as 
follows (see \cite{ye} for more detail).

\pagebreak

\bt \l{ye} Let $X$ and $P$ be non-special rational functions of degree at least two. Then $\mu_X=\mu_P$ if and only if there exist $k,l\geq 1$ such that the equalities 
\be \l{sisis} X^{\circ 2k}=X^{\circ k}\circ P^{\circ l},\ \ \ \ \ \ \ \ \ P^{\circ 2l}=P^{\circ l}\circ X^{\circ l},\ee
hold. \qed 
\et

Setting $F=X^{\circ k},$ $G=P^{\circ l},$ we can rewrite system \eqref{sisis} in the form 
\be \l{eq} F\circ F=F\circ G, \ \ \ \ G\circ G =G\circ F.\ee  The problem of describing rational solutions of this  system with $F\neq G$ is closely related to the problem of describing rational solutions of the functional equation 
\be \l{fe} A\circ X=A\circ Y, \ee distinct from the trivial solution $X=Y.$ 
Specifically, it was observed in the paper \cite{ye} that if $X$, $Y$, and $A$ are  rational functions such that equality \eqref{fe} holds, then the functions \be \l{sys} F=X\circ A, \ \ \  \ G=Y\circ A\ee
satisfy \eqref{eq}. Moreover, it was proved in \cite{entr} that {\it all} solutions of \eqref{eq} 
can be obtained in this way.  

A comprehensive classification of  rational functions satisfying \eqref{fe} is not known. The most complete result in this direction,  obtained in the paper \cite{az}, is the classification of solutions of \eqref{fe} under the assumption that $A$ is a polynomial. 
For some partial results we refer the reader to \cite{az}, \cite{entr}, \cite{r4},  \cite{seg}.  
It is instructive to consider the following more general problem. 
Let $A$ be a rational function of degree at least two. We say that $A$ is {\it tame} if the algebraic curve 
\be \l{cira} \frac{A(x)-A(y)}{x-y}=0\ee has no factors of genus zero or {\it one}. Otherwise, we say that $A$ is {\it wild}.

\bpr \l{wild} Describe wild rational functions. 
\epr

Note that by the Picard theorem, the condition that $A$ is tame is equivalent to the condition that 
 equality  \eqref{fe}, where $X$, $Y$ are functions {\it meromorphic} on $\C$,  implies that $X\equiv Y.$ 

Let us remark that Problem \ref{wild} is closely related to the question about possible decompositions of ``even'' rational functions (cf. \cite{bear}, \cite{hu1}, \cite{hu2}).  
 Let $F$ be an ``even'' rational function, that is, a rational function of the form $F=U\circ z^2,$ where $U\in \C(z)$, and 
$F=A\circ X$ an arbitrary decomposition of $F$ into a composition of rational functions.  Then the equality 
$$F=U\circ z^2=A\circ X$$ implies that  \eqref{fe} holds for
 $Y(z)=X(-z).$ Of course, if the rational function $X$ is also even, then $Y=X$. However, if $X$ is not even,
 we obtain a non-trivial solution of \eqref{fe}.   This construction can be generalized to the case where instead of $z^2$ any Galois covering  from a torus or $\C\P^1$ to $\C\P^1$  is used (see \cite{tame}). Moreover, if curve \eqref{cira} is irreducible, then all solutions of \eqref{fe}  can be obtained in this way (\cite{entr}).

\pagebreak

Rational functions sharing an iterate share a measure of maximal entropy, and the system \eqref{sisis} can be regarded as a generalization of the condition that $X$ and $P$ share an iterate. Accordingly, Conjecture \ref{co3} is a particular case of the following conjecture.

\bco \l{co5}  Let $P$ be a non-special rational function of degree at least two. Then there 
exist {\it finitely many} rational functions $X_1, X_2,\dots, X_r$ such that $X$ belongs to $E(P)$ if and only if 
\be \l{rep3} X=X_j\circ P^{\circ k}\ee
for some  $j,$ $1\leq j \leq r,$  and  $k\geq 0.$ 
\eco 

Note that Conjecture \ref{co5} implies Theorem \ref{ye} in the same way as Theorem \ref{gr0} implies the Ritt theorem. Indeed, if Conjecture \ref{co5} is true and $X\in E(P)$, then equality \eqref{yuy} holds for some $l_2>l_1\geq 1$ and $k_2>k_1\geq 1$. 
By symmetry, also 
\be \l{tog1} P^{\circ l_2'}=P^{\circ l_1'}\circ  X^{\circ (k_2'-k_1')}\ee 
for some $l_2'>l_1'\geq 1$ and $k_2'>k_1'\geq 1$. Finally, equalities \eqref{yuy} and \eqref{tog1} imply that equalities \eqref{sisis} hold 
for some $k,l\geq 1$ (see \cite{amen}, Lemma 2.10). 
Of course, in distinction with the commutative case, now we cannot exchange the order of iterates of $X$ and $P$ in equalities \eqref{yuy} or \eqref{tog1} and obtain  equality \eqref{aob}.

\section{Reversible semigroups}
It follows from the definition \eqref{reva}  that  
if a semigroup of rational functions $S$ is left reversible, then for all $A,B\in S$   the  ``separated variable" curve 
\be \l{cu0} A(x)-B(y)=0\ee
 has a factor of genus zero. Separated variable curves with a factor  of genus zero have been intensively studied  (see e. g. \cite{az}, \cite{bilu}, \cite{DLS}, \cite{f1}, \cite{f3}, 
 \cite{lau}, \cite{cur}, \cite{plo}), but their  
 full description is not known. 
Notice also that  the irreducibility problem for curves \eqref{cu0}, the so-called Davenport-Lewis-Schinzel problem,  is also very difficult and is widely open   
(see  \cite{CCN99}, \cite{DLS},  \cite{f2},   \cite{fc}, \cite{nef2}).

More generally, since $S$ is a semigroup, the left reversibility condition implies that for all $A,B\in S$  all  algebraic curves 
\be \l{cu1} A^{\circ n}(x)-B(y)=0, \ \ \ n\geq 1,\ee
 have a factor of genus zero. On the other hand, the problem of describing  $A$ and $B$ such that all curves \eqref{cu1} have a factor of genus zero or {\it one}  is a geometric counterpart of the following problem of the arithmetic nature posed in \cite{jones}:   which rational functions $A$
defined over a number field $K$ have a $K$-orbit containing infinitely many points from the value set $B(\P^1(K))$? 
These problems have been studied in the papers \cite{jones}, \cite{hyde}, \cite{aol}.  In particular, in \cite{aol}, a description of such $A$ and $B$  in terms of semiconjugacies and Galois coverings was obtained. Notice that all curves \eqref{cu1} obviously have a factor of genus zero whenever $B$ is a ``compositional left factor'' of some iterate of $A$, where by  a compositional left factor of a rational function $F$ we mean 
any rational function  $G$ such that $F=G\circ H$ for some rational function $H.$ Moreover, 
in case $A$ is non-special, the main result of \cite{aol} in a slightly simplified form  can be formulated as follows (see \cite{aol}, Theorem 1.2).

\bt \l{tt1} Let  $A$ be a  non-special rational function of degree at least two.
 Then  there exist rational functions  $X$ and  $F$ such that $X$ is  a Galois covering,  
the diagram 
\be \l{udod}
\begin{CD}
\C\P^1 @>F>> \C\P^1 \\
@VV X V @VV  X V\\ 
\f \C\P^1 @>A>> \f\C\P^1\
\end{CD}
\ee
commutes, and for a rational function $B$ of degree at least two all algebraic curves \eqref{cu1}  have a factor of genus zero or one   if and only if $B$ is a compositional left factor  of $A^{\circ \ell}\circ X $ for some $l\geq 0.$ 
\et 

Finally, we observe that if a semigroup of rational functions $S$ is left reversible, then the following condition holds: 
for all $A,B\in S$ all  algebraic curves 
\be \l{cu2} A^{\circ n}(x)-B^{\circ m}(y)=0, \ \ \ n,m\geq 1,\ee
 have a factor of genus zero. This condition is stronger than the previous two conditions, and we conjecture that for such $A$ and $B$ the following stronger version of Theorem \ref{tt1} holds. 

\bco \l{41} Let  $A$ and $B$ be rational functions of degree at least two such that  all algebraic curves \eqref{cu2}  have a factor of genus zero or one. Then either both $A$ and $B$ are special or 
there exist $k,l\geq 1$ such that $A^{\circ k}=B^{\circ l}$.
\eco

If true, Conjecture \ref{41} implies  the following conjecture. 
\bco \l{42} Let $S$ be a semigroup of rational functions of degree at least two containing at least one non-special function.  Then $S$ is left reversible if and only if any two elements of $S$ have a common iterate.

\eco

Indeed, Conjecture \ref{41} implies  the ``only if'' part of Conjecture \ref{42}. On the other hand, the ``if'' part is obvious, since the equality $A^{\circ k}=B^{\circ l}$ implies that \eqref{reva} is satisfied for $X=A^{\circ (k-1)}$ and $Y=B^{\circ (l-1)}.$ Moreover, since $A^{\circ k}=B^{\circ l}$ implies 
$A^{\circ 2k}=B^{\circ 2l}$, we can assume that $k,l\geq 2$ ensuring  that $X\in S$ and $Y\in S.$

Furthermore, Conjecture \ref{41} implies the following conjecture.

\bco  \l{43} Let $A$ and $B$ be  rational functions of degree at least two such that an orbit of $A$ has an infinite intersection with an 
 orbit of $B$. Then $A$ and $B$ have a common iterate. 
\eco

 In case $A$ and $B$ are polynomials, Conjecture \ref{42} is the theorem proved in the papers \cite{gtz}, \cite{gtz2}. This result was extended to tame rational functions in the paper \cite{tame}.
Similarly, Conjecture \ref{41} and Conjecture \ref{42} are true if the functions involved are polynomials or tame rational functions 
(\cite{amen}).

To see that Conjecture \ref{41} implies Conjecture \ref{43}, we recall that, by the Faltings theorem (\cite{fa}), if an irreducible algebraic curve $C$ defined over a finitely generated field $K$ of characteristic zero has infinitely many $K$-points, 
then $g(C)\leq 1.$  On the other hand, it is easy to see that if the orbit intersection $O_A(z_1)\cap O_B(z_2)$ is infinite,  then all curves \eqref{cu2} have infinitely many points $(x,y)\in O_A(z_1)\times O_B(z_2).$ Since the orbits $O_A(z_1)$, $ O_B(z_2)$ belong to the field $K$ finitely generated over $\Q$ by $z_1$, $z_2$, and the coefficients of $A$, $B$, this implies that   all  curves \eqref{cu2} have  
a factor of genus zero or one. Taking into account that Conjecture \ref{43} is true for special $A$ and $B$ (\cite{tame}), this shows that 
Conjecture \ref{41} implies Conjecture \ref{43}.

Switching to right reversible semigroups, instead of the condition that 
for all $A,B\in S$   the algebraic curve \eqref{cu0} has a factor of genus zero, we obtain the condition that for all $A,B\in S$   the field 
$\C(A)\cap \C(B)$ contains a non-constant rational function. Thus, 
 we face the following  problem.

\bpr Given rational functions $A$ and $B$ of degree at least two, under what conditions  does the field 
$\C(A)\cap \C(B)$ contain a non-constant rational function?
\epr

Despite a very natural setting of this problem, essentially nothing is known about its solutions unless both $A$ and $B$ are polynomials, in which case a complete description of such $A$ and $B$ is known. Specifically, if $\C(A)\cap \C(B)$ contains a {\it polynomial} the answer is given by the Ritt theory (\cite{r1}), and the general case reduces to this  one (\cite{pcon}). 
For some other related results we refer the reader to the papers \cite{ac}, \cite{bt}, \cite{hs}, \cite{lys}.

The analogues of the other two problems about algebraic curves considered above 
can be formulated as follows: given rational functions $A$ and $B$, under what conditions all the fields  \be \l{fi1} \C(A^{\circ n})\cap \C(B), \ \ \ n\geq 1,\ee and, more generally, all the fields  
\be \l{fi2} \C(A^{\circ n})\cap \C(B^{\circ m}), \ \ \ n,m\geq 1,\ee do contain a non-constant rational function?

We conjecture that the first problem has the following solution,  ``symmetric'' to the one provided by Theorem \ref{tt1}.  

\bco \l{coo} Let  $A$ be a  non-special rational function of degree at least two.
 Then  there exist rational functions  $X$ and  $F$ such that  $X$ is  a Galois covering,  
the diagram 
\be \l{xor} 
\begin{CD}
\C\P^1 @>A>> \C\P^1 \\
@VV X V @VV  X V\\ 
\f \C\P^1 @>F>> \f\C\P^1\
\end{CD}
\ee
commutes, and for a rational function $B$ of degree at least two all fields \eqref{fi1} contain a  non-constant rational function   if and only if $B$ is a compositional right factor  of $ X\circ A^{\circ \ell} $ for some $l\geq 0.$

\eco

Notice that the ``if'' part of  Conjecture \ref{coo} is obtained by a direct calculation. 
Indeed, it follows from \eqref{xor} and   $$ X\circ A^{\circ l}=U\circ B$$ 
that for every $k\geq 0$ the equality 
$$F^{\circ k}\circ U\circ B=F^{\circ k}\circ  X\circ A^{\circ l}=X\circ A^{\circ l+k}$$ holds, implying that 
 all fields \eqref{fi1} contain a  non-constant rational function.  
In particular,  they contain a  non-constant rational function whenever $U$  is
a  compositional right factor of some iterate  $A^{\circ l}$, $l\geq 1$ (the case where $F=A$ and $X=z$).

Finally, analogues of Conjecture \ref{41} and Conjecture \ref{42} are the following conjectures, which are known to be true in the polynomial case (\cite{amen}).

\bco  Let  $A$ and $B$ be rational functions of degree at least two such that   all fields \eqref{fi1} 
contain a non-constant rational function.  Then either both $A$ and $B$ are special or 
there exist $k,l\geq 1$  such that the equalities $$A^{\circ 2k}=A^{\circ k}\circ B^{\circ l}, \ \ \  \ B^{\circ 2l}=B^{\circ l}\circ A^{\circ k}$$ hold.  
\eco

\bco  Let $S$ be a semigroup of rational functions of degree at least two containing at least one non-special function.  Then $S$ is right reversible if and only if for any two elements $A,B$ of $S$  there exist $k,l\geq 1$  such that the equalities $$A^{\circ 2k}=A^{\circ k}\circ B^{\circ l}, \ \ \  \ B^{\circ 2l}=B^{\circ l}\circ A^{\circ k}$$ hold.

\eco

\section{Amenable semigroups}

The following conjecture, generalizing the corresponding result for polynomials proved in \cite{amen}, presumably describes amenable and left amenable semigroups of rational functions.

\bco \l{ass1} 

Let $S$ be a semigroup of rational functions of degree at least two containing at least one non-special rational function.  Then the following conditions are equivalent.

\begin{enumerate} [label=\arabic*)]

\item The semigroup $S$   is   amenable.

\item The semigroup $S$   is   left amenable.

\item The semigroup $S$   is    left reversible.

\item Any two elements of $S$  have a common iterate.

\item  The semigroup $S$   is  a subsemigroup of $C_{\infty}(P)$  
  for some non-special rational function $P$ of degree at least two.   	
\end{enumerate} 

\eco

Notice that to prove Conjecture \ref{ass1} it is enough to prove  the implication $3\Rightarrow 4$ only. Indeed, the implication $1\Rightarrow 2$ is obvious, the implication $2\Rightarrow 3$ is well-known, and it is easy to see 
that if $4)$ holds, then $S\subseteq C_{\infty}(P)$ for any  $P\in S$, implying the implication $4\Rightarrow 5$. Finally, the implication $5\Rightarrow 1$ is proved in \cite{amen}. 
Note that since Conjecture \ref{41} implies through Conjecture \ref{42} the implication $3\Rightarrow 4,$ to prove Conjecture \ref{ass1} it is enough to prove Conjecture \ref{41}.

The next conjecture describes right amenable semigroups of rational functions, generalizing results proved in \cite{peter}, \cite{amen} in the  polynomial case.  

\bco \l{ass2} 

Let $S$ be a semigroup of rational functions of degree at least two containing at least one non-special rational function.  Then the following conditions are equivalent.

\begin{enumerate} [label=\arabic*)]

\item The semigroup $S$   is right amenable.

\item The semigroup $S$   is   right reversible.

\item  For any two elements $A,B$ of $S$  there exist $k,l\geq 1$  such that the equalities $$A^{\circ 2k}=A^{\circ k}\circ B^{\circ l}, \ \ \  \ B^{\circ 2l}=B^{\circ l}\circ A^{\circ k}$$ hold.

\item The semigroup $S$ contains no free subsemigroup of rank two.

\item The semigroup $S$ is a subsemigroup of  $E(P)$ for some  non-special rational function $P$   
of degree at least two.

\end{enumerate} 

\eco

Let us list the implications in Conjecture \ref{ass2}, which are known to be true. 
The implication $1\Rightarrow 2$ is well-known. The implications $3\Rightarrow 4$ and $3\Rightarrow 2$ are obvious. The equivalence $3\Leftrightarrow 5$ follows easily from Theorem \ref{ye}. Finally,   the implication $4\Rightarrow 2$ can be established as follows (cf. \cite{amen}, Lemma 2.8). 
Let $A\neq B$ be elements of $S$. By condition, the semigroup generated by $A$ and $B$ is not free. Therefore, there exist $X_1,X_2,\dots ,X_k\in \{A,B\}$ and   
$Y_1,Y_2,\dots ,Y_l\in \{A,B\}$ such that \be \l{sok} X_kX_{k-1}\dots X_1=Y_lY_{l-1}\dots Y_1,\ee but the words in $A$ and $B$ in the parts of this equality are different. 
It follows  from the right cancellativity of $\C(z)$ that without loss of generality we may assume that $X_1\neq Y_1$. Thus, either $X_1=A$ and $Y_1=B$, or $X_1=B$ and $Y_1=A.$ Moreover, the condition $A\neq B$ implies that $k\geq 2$ and $l\geq 2$. Thus, if, say, $X_1=A$ and $Y_1=B$, the elements  
$$X= X_k\dots X_2, \ \ \ \ \   Y=Y_l\dots Y_2 $$ of $S$ satisfy $X\circ A=Y\circ B$.

\end{document}